\documentclass{svmult}

\usepackage{makeidx}         
\usepackage{graphicx}        
\usepackage{multicol}        
\usepackage[bottom]{footmisc}
\usepackage{natbib}
\usepackage{url}
\usepackage{xcolor}
\usepackage{subfig}

\renewcommand{\PackageWarningNoLine}[2]{}
\newcommand{\tvect}[1]{{\bf #1}}

\newcommand{\dt}{\Delta t}
\renewcommand{\d}{\mathrm{d}}

\usepackage{amsmath}
\usepackage{amsfonts}
\usepackage{array}

\begin{document}

\title*{Inexact spectral deferred corrections}
\author{Robert Speck\inst{1,2}\and
Daniel Ruprecht\inst{2}\and
Michael Minion\inst{3}\and
Matthew Emmett\inst{4}\and
Rolf Krause\inst{2}}
\authorrunning{R.~Speck et al.}
\institute{J\"ulich Supercomputing Centre, Forschungszentrum J\"ulich GmbH, Germany
\texttt{r.speck@fz-juelich.de}\and
Institute of Computational Science, Universit{\`a} della Svizzera italiana, Lugano, Switzerland
\texttt{daniel.ruprecht@usi.ch, rolf.krause@usi.ch}\and
Institute for Computational and Mathematical Engineering, Stanford University,
Stanford, USA \texttt{mlminion@gmail.com}\and
Center for Computational Sciences and Engineering, Lawrence Berkeley National Laboratory, USA
\texttt{mwemmett@lbl.gov}}
%
%
\maketitle

\section{Introduction}\label{speck_mini_15_sec:intro}
Implicit integration methods based on collocation are attractive for a number of reasons, e.g. their ideal (for Gauss-Legendre nodes) or near ideal (Gauss-Radau or Gauss-Lobatto nodes) order and stability properties.
However, straightforward application of a collocation formula with $M$ nodes to an initial value problem with dimension $d$ requires the solution of one large $Md \times Md$ system of nonlinear equations.

Spectral deferred correction (SDC) methods, introduced by~\cite{DuttEtAl2000}, are an attractive approach for iteratively computing the solution to the collocation problem using a low-order method (like implicit or IMEX Euler) as a building block. Instead of solving one huge system of size $Md \times Md$, SDC iteratively solves $M$ smaller $d \times d$ systems to approximate the solution of the full system (see also the discussion in~\cite{HuangEtAl2006}).
It has been shown e.g. by~\cite{ShuEtAl2007} that each iteration/sweep of SDC raises the order by one, so that SDC with $k$ iterations and a first-order base method is of order $k$, up to the order of the underlying collocation formula.
Therefore, to achieve formal order $p$, SDC requires $p/2$ nodes and $p$ iterations and thus $p^{2}/2$ solves of a $d \times d$ system (for Gauss-Legendre nodes).

Considering the number of solves required to achieve a certain order, one might conclude that, notwithstanding the results presented here, SDC is less efficient than e.g. diagonally implicit Runge Kutta (DIRK) methods, see e.g. \cite{Alexander1977}, which only require $p-1$ solves.
However, the flexibility of the choice of the base propagator in SDC and the very favorable stability properties make it an attractive method nevertheless.
In particular, semi-implicit methods of high order can easily be constructed with SDC which make it competitive for complex applications, see~\cite{Minion2003} and~\cite{BourliouxEtAl2003}.
Further extensions to SDC allow it to integrate processes with different time scales, see~\cite{LaytonMinion2004,bouzarthMinion:2011} efficiently; and the iterative nature of SDC also allows it to be extended to a multigrid-like multi-level algorithm, where work is shifted to coarser, computationally cheaper levels, see~\cite{SpeckEtAl2013_mlsdc}.

In the present paper, we introduce another strategy to improve the efficiency of SDC, which is
similar to ideas from~\cite{oosterlee1996use} where a single V-cycle of a multigrid method is used as a preconditioner.
We show here that the iterative nature of SDC allows us to use incomplete solves of the linear systems arising in each sweep.
In the resulting \emph{inexact spectral deferred corrections} (ISDC), the linear problem in each Euler step is solved only approximately using a small number (two in the examples presented here) of multigrid V-cycles.
It is numerically shown that this strategy results in only a small increase of the number of required sweeps while reducing the cost for each sweep.
We demonstrate that ISDC can provide a significant reduction of the overall number of multigrid V-cycles required to complete an SDC time step.

\section{Semi-implicit spectral deferred corrections}
We consider an initial value problem in Picard form
\begin{align}
   u(t) = u_0 + \int_{T_0}^{t}f(u(s))\ \d s
\end{align}
where $t\in[T_0,T]$ and $u,f(u)\in\mathbb{R}^N$.
Subdividing a time interval $[T_n,T_{n+1}]$ into $M$ intermediate substeps $T_n = t_0 \le t_1<...<t_M\le T_{n+1}$, the integrals from $t_m$ to $t_{m+1}$ can be approximated by
\begin{align}
   I_m^{m+1} = \int_{t_m}^{t_{m+1}}f(u(s))\ \d s \approx \dt\sum_{j=0}^Ms_{m,j}f(u_j) = S_m^{m+1}F(u)
\end{align}
where $u_m \approx u(t_m)$, $m=0,1,...,M$, $F(u) = (f(u_1),...,f(u_M))^T$, $\dt = T_{n+1}-T_n$, and $s_{m,j}$ are quadrature weights.
The nodes $t_m$ correspond to quadrature nodes of a spectral collocation rule like Gauss-Legendre or Gauss-Lobatto quadrature rule.
The basic implicit SDC update formula at node $m+1$ in iteration $k+1$ can be written as
\begin{align}
   u_{m+1}^{k+1} = u_m^{k+1} + \dt_m\left[f(u_{m+1}^{k+1})-f(u_{m+1}^k)\right] + S_m^{m+1}F(u^k),
\end{align}
where $\dt_m = t_{m+1}-t_m$, for $m=0,...,M-1$.
Alternatively, if $f$ can be split into a stiff part $f^I$ and a non-stiff part $f^E$, a semi-implicit update is easily
constructed for SDC using
\begin{align}
   u_{m+1}^{k+1} = u_m^{k+1} &+ \dt_m\left[f^I(u_{m+1}^{k+1})-f^I(u_{m+1}^k)\right]\nonumber \\
   &+ \dt_m\left[f^E(u_{m}^{k+1})-f^E(u_{m}^k)\right] + S_m^{m+1}F(u^k)\label{speck_mini_15_eq:sisdc}.
\end{align}
Here, only the $f^I$-part is treated implicitly, while $f^E$ is explicit.
We refer to \cite{Minion2003} for the details on semi-implicit spectral deferred corrections.

\section{Inexact spectral deferred corrections}
In the following, we consider the linearly implicit case $f^I(u) = Au$, where $A$ is a discretization of the Laplacian operator.
Here, spatial multigrid is a natural choice for solving the implicit part in~\eqref{speck_mini_15_eq:sisdc}.
As in~\cite{SpeckEtAl2013_mlsdc}, we use a high-order compact finite difference stencil to discretize the Laplacian
(see e.g. ~\cite{spotzcarey_1996}).
This results in a weighting matrix $W$ for the right-hand side of the implicit system and, with the notation $\tilde{f}^I(u) = Wf^I(u)$, the semi-implicit SDC update~\eqref{speck_mini_15_eq:sisdc} becomes
\begin{align}
   (W-\dt_mA)u_{m+1}^{k+1} = Wu_{m}^{k+1} &+ \dt_mW\left[f^E(u_m^{k+1})-f^E(u_m^{k})\right]\nonumber\\ &- \dt_m\tilde{f}^I(u_{m+1}^k) + S_m^{m+1}\tilde{F}(u^k),
\end{align}
where $\tilde{f} = Wf^E + \tilde{f}^I$ and  $\tilde{F}(u^k) = (\tilde{f}(u_1^k),...,\tilde{f}(u_M^k))^T$.
Thus, instead of inverting the operator $I-\dt_mA$ in~\eqref{speck_mini_15_eq:sisdc}, the right-hand side of~\eqref{speck_mini_15_eq:sisdc} is modified by $W$ and the operator $W-\dt_mA$ needs to be inverted.
We note that for calculating the residual during the SDC iteration, the weighting matrix needs to be inverted once per node, which can be done using multigrid as well.

For classical SDC, each computation of $u_{m+1}^{k+1}$ includes a full inversion of $W-\dt_mA$ using e.g.~a multigrid solver in space.
For $K$ iterations and $M$ nodes, the multigrid solver is executed $K(M-1)$ times, each time until a predefined tolerance is reached.
In order to reduce the overall number of required multigrid V-cycles, ISDC replaces this full solve with a a small fixed number $L$ of V-cycles, leading to an accumulated number of $\tilde{K}(M-1)L$ V-cycles in total.
Naturally, the number of iterations in ISDC will be larger than the number of SDC, that is $K\le\tilde{K}$.
However, if $\tilde{K}$ is small enough so that $\tilde{K}(M-1)L$ is below the total number of multigrid V-cycles required for
$K(M-1)$ full multigrid solves, inexact SDC will be more efficient than classical SDC.

Convergence is monitored using the maximum norm of the SDC residual, a discrete analogue of $u^{k}(t) - u_0 - \int_{T_0}^{t} f(u^{k}(s))\ \d s$, that measures how well our iterative solution satisfies the discrete collocation problem.
See~\cite{SpeckEtAl2013_mlsdc} for definition and details.
In the tests below, sweeps are performed until the SDC residual is below a set threshold.

\section{Numerical tests}
In order to illustrate the performance of ISDC, we consider two different numerical examples, the 2D diffusion equation and 2D viscous Burgers' equation.
As described above, in both cases the diffusion term is discretized using a $4$th-order compact stencil with weighting matrix and a spatial mesh with $64$ points.
For Burgers' equation, the advection term is is discretized using a $5^{\rm th}$ order WENO scheme.
\subsection{Setups}
The first test problem is the 2D heat equation on the unit square, namely
\begin{align}
   u_t(\tvect{x},t) &= \nu \Delta u(\tvect{x},t),\quad \tvect{x} \in \Omega = [0,1]^2\\
   u(\tvect{x},0) &= \sin(\pi x) \sin(\pi y)\\
   u(\tvect{x},t) &= 0 \ \mbox{on} \ \partial \Omega
\end{align}
with $\tvect{x} = (x,y)$.
The exact solution is $u(x,t) = \exp(-2 \pi^2 \nu t) \sin(\pi x) \sin(\pi y)$.
An implicit Euler is used here as base method in SDC.

The second test problem is the nonlinear viscous Burgers' equation
\begin{align}
   u_t(\tvect{x},t) + u(\tvect{x},t) u_{x}(\tvect{x},t) + u(\tvect{x},t) u_{y}(\tvect{x},t)  &= \nu \Delta u(\tvect{x},t), \tvect{x}\in[-1,1]^2, \\
   u(x,0) &= \exp\left(-\frac{x^2}{\sigma^2}\right),\quad \sigma = 0.1
\end{align}
with periodic boundary conditions.
Here, an IMEX Euler is used as base method, i.e. the Laplacian is treated implicitly, while the advection term is integrated explicitly.

In both examples the diffusion parameter $\nu$ controls the stiffness of the term $f^I$: for a given spatial resolution, the shifted Laplacian $W-\nu\dt A$, and therefore the performance of the multigrid solver, depends critically on $\nu$.
We choose three different values of $\nu$ for each example to measure the impact of stiffness on the performance of ISDC: $\nu=1$, $10$, $100$ for the heat equation and $\nu = 0.1$, $1$, $10$ for Burgers' equation.
For ISDC, each implicit solve is approximated using $L=2$ V-cycles.
A single time-step of length $\Delta t = 0.001$ is analyzed for a spatial discretization with $\Delta x = \Delta y = 1/64$ in both cases, leading to CFL numbers for the diffusive term of approximately $4.1$, $41$ and $410$ for the heat equation and $0.41$, $4.1$ and $41$ for Burgers' equation.

\subsection{Results}
\begin{table}[t]
\centering
\subfloat[Heat equation]{
\begin{tabular}{c|r @{ } r @{ \ } r @{ \ } rr>{\bf}r}
 $\nu$         & & $M$  &  SDC & ISDC &  & savings      \\ \hline
                   &  &           &           &            &  &        \\
$1$             &  &  3       & 16(4)  &  12(4) &  & 25\%  \\
                   &  &  5       & 23(3)  &  20(3) &  & 13\% \\
                   &  &  7       & 32(3)  &  28(3) &  & 13\% \\
                   &  &           &           &            &  &        \\
$10$           &  & 3        & 36(5)  &  20(5) &  & 44\% \\
                   & &  5        & 61(5)  & 40(5)  &  & 34\% \\
                   & &  7        & 79(4)&   47(4)  &  & 41\% \\
                   & &            &           &            &  &         \\
  $100$       & & 3        & 106(13) & 52(13)   & & 51\% \\
                   & & 5        & 150(10) & 104(13) & & 31\% \\
                   & & 7        & 187(9)\phantom{0}  & 167(14) & & 11\% \\
                   & &           &           &             & &         \\
\end{tabular}
}\qquad\qquad
\subfloat[Viscous Burgers' equation]{
\begin{tabular}{c|r @{ } r @{ \ } r @{ \ } rr>{\bf}r}
 $\nu$        & & $M$  &  SDC & ISDC & & savings          \\ \hline
                   &  &           &                &             &  &        \\
$10^{-1}$    &  &   3      & 21(8)      & 21(8)   &   & 0\%  \\
                   &  &   5       & 26(6)     & 26(6)   &   &  0\% \\
                   &  &   7       & 33(5)     & 33(5)   &   &  0\% \\
                   &  &            &              &             &   &        \\
$1.0$          &  &   3       & 97(17)   & 66(17) &   & 32\% \\
                   & &    5       & 140(17) & 117(17) &  & 16\% \\
                   & &    7       & 160(15) & 143(15) &  & 11\% \\
                   & &            &               &               &  &         \\
  $10$        & & 3         & 207(25)   & 100(25)  & & 52\% \\
                   & & 5        & 523(38)   & 298(38)  & & 43\% \\
                   & & 7        & 902(50)   & 578(50)  & & 36\% \\
                   & &           &           &             & &         \\
\end{tabular}
}
\caption{Accumulated multigrid V-cycles for (a) the heat equation and (b) the viscous Burgers' equation with different values for the diffusion coefficients $\nu$ and the number  of quadrature nodes $M$. Cycles are accumulated over all sweeps required to reduce the SDC or ISDC residual below $5 \times 10^{-8}$. The number of deferred correction sweeps is shown in parentheses.
Saving indicates the amount of V-cycles saved by ISDC in percent of the cycles required by SDC.}
\label{speck_mini_15_tab:ncycles}
\end{table}
Table~\ref{speck_mini_15_tab:ncycles} shows the total number of multigrid V-cycles for the heat equation (left) and for Burgers' equation (right) for three different numbers of collocation nodes $M$ and different values of $\nu$.
The number of SDC or ISDC sweeps performed is shown in parentheses.
In each case, sweeps are performed until the SDC or ISDC residual is below $5 \times 10^{-8}$.
To simplify the analysis in the presence of the weighting matrix, the V-cycles required to invert the weighting matrix are not counted here.
In the last row, the amount of V-cycles saved by ISDC is given in percent of the required SDC cycles.
%

In most cases,  ISDC provides a substantial reduction of the total number of required multigrid V-cycles and requires only slightly more sweeps to converge than SDC.
The most savings can be obtained if the number of multigrid V-cycles in SDC is high but ISDC does not lead to a significant increase in sweeps, which is the case for mildly stiff problems (e.g. $\nu=10$ for heat equation and $\nu=1$ and $\nu=10$ for Burgers) or stiff problems with small values for $M$.
For stiff problems with large $M$ (e.g. heat equation with $\nu=100$ and $M=7$), however, ISDC leads to a more significant increase in required sweeps, therefore only resulting in small savings.
For the non-stiff cases, particularly for Burger's equation, ISDC does not provide much benefit, but also does no harm: the multigrid solves in SDC take only very few V-cycles to converge, so that SDC and ISDC are almost identical (for Burgers with $\nu=0.1$, SDC and ISDC are actually identical).
In a sense, for simple problems where the stopping criterion of the multigrid solver is reached after one or two V-cycle anyhow, SDC automatically reduces to ISDC.

In summary, the tests presented here suggest that replacing full multigrid solves by a small number of V-cycles in SDC only leads to a small increase in the total number of SDC sweeps required for convergence but can significantly reduce the
computational cost of each sweep. The savings in the overall number of multigrid V-cycles of ISDC directly translates into faster run times of ISDC runs compared to classical SDC.

Preliminary numerical tests not document here suggest that, as long as the approximate solution of the linear system is sufficiently accurate, the order of ISDC still increases by each iteration, as shown for SDC in~\cite{ShuEtAl2007}. 
A detailed study confirming this, including a possible extension of the proof, is left for future work.
\subsection{Interpretation}
The good performance of ISDC in the examples presented above is mainly due to the choice of the starting values for the multigrid solver.
When performing the implicit Euler step to compute $u_{m+1}^{k+1}$, the value $u^{k}_{m+1}$ from the previous SDC sweep gives a very good starting value, particularly in later sweeps.
Therefore, even two multigrid V-cycles are sufficient to approximate the real solution of the linear system of equations reasonably well. This effect can be observed by monitoring the number of V-cycles in classical SDC. During the first sweep, many more V-cycles
are typically required for multigrid to converge than in later sweeps where the initial guess becomes  very accurate as the SDC
iterations converge.  In fact, during the last sweeps of SDC, a single V-cycle is often sufficient for solving the implicit system.
Hence, the additional sweeps required by ISDC are  mainly due to the less accurate approximations during the first sweeps.
As soon as the initial guess $u_{m+1}^{k}$ for $u_{m+1}^{k+1}$ is good enough, ISDC basically proceeds like SDC.
A computational experiment that confirms this is as follows: if, when solving for $u_{m+1}^{k+1}$, we replace the initial guess with the zero vector, or even $u_{m}^{k+1}$, then
ISDC fails to converge altogether.
On the other hand, SDC still convergences in this scenario, but the number of required multigrid V-cycles increases dramatically.

It is important to  contrast this behavior to non-iterative schemes like diagonally implicit Runge-Kutta, where usually only the value from the previous time step or stage is available to be used as starting value.  Our experience with SDC methods suggests that more
multigrid V-cycles would be required to solve each stage in  a DIRK scheme than in later SDC iterations.  Hence, simply
counting the number of implicit function evaluations required could be a misleading way to compare the cost of SDC
and DIRK schemes.

\section{Conclusion and Outlook}
The paper presents a variant of spectral deferred corrections called \emph{inexact spectral deferred corrections} that can significantly reduce the computational cost of SDC.
In ISDC, full spatial solves within SDC sweeps with an implicit or semi-implicit Euler are replaced by only a few V-cycles of a multigrid.
In the two investigated examples, ISDC saves up to $52$\% of the total multigrid V-cycles required by SDC with full linear solves in each step, while only minimally increasing the number of sweeps required to reduce the SDC residual below some set tolerance.
The main reason for the good performance of ISDC is that the iterative nature of SDC provides very accurate initial guesses for the multigrid solver.
Besides providing significant speedup, ISDC essentially removes the need to define a tolerance or maximum number of iterations for the spatial solver.

A natural extension of the work presented in this paper is the application of ISDC sweeps in MLSDC, the multi-level version of SDC.
MLSDC performs SDC sweeps in  a multigrid-like way on multiple levels.
The levels are connected through an FAS correction term in forming the coarsened  spatial representation of the problem on upper levels of the hierarchy.
Using ISDC corresponds to the ``reduced implicit solve'' strategy mentioned in~\cite{SpeckEtAl2013_mlsdc} and incorporating it into MLSDC could further improve its performance.
Finally, the ``parallel full approximation scheme in space and time'' (PFASST, see~\cite{Minion2010,EmmettMinion2012,EmmettMinion2012_DDM} for details) performs SDC sweeps on multiple levels combined with a forward transfer of updated initial values in a manner similar to Parareal (see~\cite{LionsEtAl2001}).
Instead of performing a full time integration as done in Parareal, PFASST interweaves SDC sweeps with Parareal iterations so that on each time level, only a single SDC sweep is performed (i.e.~an inexact time integrator is applied), leading to a time-parallel method with good parallel efficiency (see e.g.~\cite{SpeckEtAl2013_Parco} and \cite{RuprechtEtAl2013_SC}).
Integrating ISDC into PFASST could further improve its parallel efficiency.

\subsubsection*{Acknowledgments}
This work was supported by the Swiss National Science Foundation (SNSF) under the lead agency agreement through the project ``ExaSolvers'' within the Priority Programme 1648 ``Software for Exascale Computing'' (SPPEXA) of the Deutsche Forschungsgemeinschaft (DFG). Matthew Emmett and Michael Minion were  supported by the Applied Mathematics Program of the DOE Office of Advanced Scientific Computing Research under the U.S. Department of Energy under contract DE-AC02-05CH11231.  Michael Minion was also supported by
the U.S. National Science Foundation grant DMS-1217080. The authors acknowledge support from Matthias Bolten, who provided the employed multigrid solver.

%
\bibliographystyle{plainnat}
\bibliography{speck_mini_15} 



\end{document}